# Lattice of Triangulations: the proof and an algorithm

Thinh Duc Nguyen [1], Duong Thi Ha Phan [2]

## Abstract

In this paper, we prove that the set of triangulations of a polygon can be equipped with an order to become a lattice. First, we define this order. In [HN99], authors defined the flip operator and then prove some properties of the graph of triangulations. We use their theorems and extend them to construct the lattice of triangulations. We prove this lattice property and introduce an elegant algorithm which correctness is induced from the proof. The complexity of this algorithm will be considered. This algorithm is efficient to find the infimum of a pair of triangulations.

*Keywords:* Triangulation; Lattice; Order; Infimum; Algorithm; Complexity

## 1. Some definitions

Consider an n-polygon *P*, we label its vertices (counter-clockwise) as *1* to *n*. In [HN99], the authors defined the flip operator that turns a triangulation into another, then they defined the graph of triangulations $G_T$. This graph is undirected. To give the order to the set of triangulations, we define the positive flip operator.

*Orders in $L_T$:*

Consider a diagonal *ac* in a triangulation $P_1$ which is inside a quadrilateral *abcd* that has ($a<b<c<d$), positive flip can be applied to diagonals of this type. In this case, it turns $P_1$ into $P_2$ which shares the same set of diagonals excetpt that *ac* is replaced by *bd*. Then we denote by $P_1$->$P_2$ to say that $P_2$ can be obtained from $P_1$ directly by a positive flip operator. This operator can be considerd as the order of triangulations.

This order creates a lattice of triangulations which we denote by $L_T(n)$. This claim will be proved in the next section.

## 2. Proof of lattice property

In [HN99], authors proved that at n*th* level of the hierarchy of triangulations, there are exactly $C_{n-2}$ vertices and this level corredponds to the $G_T$ of an n-polygon. More importantly, they showed that n*th* level contains $C_{n-3}$ paths: each vertex of the (n-1)*th* has its sons and these sons are in a path. With the above order, these paths have some more nice properties.

*Theorem 1*

Let *T* be a triangulation in *T(n)*, we have:

Each triangulation 's sons induce a subgraph on $L_T(n+1)$ that is a *directed* path which goes from left to right (the sons are sorted by their *(in)* diagonal on which the *Son* operator was applied: $S^1(T)$, $S^i(T)$,$S^j(T)$,...,$S^{n-1}(T)$ , where $1<i<j<n-1$).

[1] College of Technology, Vietnam National University, 144 Xuan Thuy-Cau Giay-Ha Noi-Vietnam
[2] Institute of Mathematics, Vietnamese Academy of Science and Technology, 18 Hoang Quoc Viet–Cau Giay–Ha Noi-Vietnam

*Proof*

Assume the vertex n is connected with some vertices (by diagonals, edges): 1, *... ,i,j, ... ,c,n-1*. The $son^j(T)$ can be obtained from $son^i(T)$ by positively flipping the edge *(n,i)* to the edge *(n+1,j)*.

Immediately, the set of directed diagonals of $L_T$ is split into 2 separate sets:

1. The above directed paths, we called a path as a family because it's induced by the sons of a common father in the directly higher level.
2. Other edges must be of the following forms: $son^i(P_1)$ -> $son^i(P_2)$ where $P_1$->$P_2$ (this property is deduced from *Lemma 3.1* in [HN99])

It means the $son^i(.)$ can only be connected to other $son^i(.)$ by a positive flip, of course except the previous and next sibling in the same family which is the path in 1. ).

Now, we will describe the proof of the conjecture. It's based on an elegant algorithm.

*Theorem 2*

$L_T(n)$ is a lattice, for every $n>=3$.

*Proof:*

It uses the inductive construction:

- Basically, when $n=3$ there is only one triangulation, so that $L_T(n)$ is a lattice.
- If $L_T(n)$ is a lattice then $L_T(n+1)$ is also a lattice. We will prove this right now.

For every pair of vertices $son^i{}_1(P_1)$ and $son^i{}_2(P_2)$ in the (n+1)*th* level, the algorithm finds their infimum $son^i{}_3(P_3)$.

$P_1$ and $P_2$ are in the n*th* level so by the induction hypothesis, they have an infimum called $P_3$.

The ways $P_1$ goes to $P_3$ and $P_2$ goes to $P_3$ is parallel to the ways their sons come to $son^i{}_1(P_3)$ and $son^i{}_2(P_3)$. (Note that the ($i_1n$) and ($i_2n$) diagonals of $P_1$ and $P_2$ cannot be positively flipped). These vertices are in one directed path which has the father $P_3$. We denote $i_3=min(i_1, i_2)$.

So, $son^i{}_3(P_3)$ is the infimum of the above 2 vertices. It implies that $L_T(n+1)$ is also a lattice.

By induction, we can conclude $L_T(n)$ is lattice for every $n>=3$. *Q.E.D*

This proof induces an algorithm for finding the infimum of a pair of triangulations in the lattice. We will describe this algorithm in the next section.

## 3. The algorithm

*Input:* A pair of triangulations $T_1$, $T_2$

*Output:* Their infimum $T_3$

*Algorithm:*

*TNDP*($T_1$, $T_2$):
If $T_1 = T_2$ return $T_1$;
Compute $P_1$, $i_1$ such that $T_1 = son^{i}_{1}(P_1)$
Compute $P_2$, $i_2$ such that $T_2 = son^{i}_{2}(P_2)$
$P_3 = TNDP(P_1, P_2)$
$i_3 = max(i_1, i_2)$
return $son^{i}_{3}(P_3)$

4. Complexity of TNDP algorithm

4.1. Trivial implementation:

If we implement the above algorithm by a trivial way, the time complexity will be $O(n^2)$. In this way, we use the *list* of diagonals and compute the $son^i(.)$ and father(.) function by trivially updating this list. We have to go through at most (n-2) levels. In each level, we compute the two function $son^i(.)$ and *father(.)* in time O(n). So, the overall time complexity is $O(n^2)$. Now, we describe the TNDP algorithm which runs in linear time.

*4.2. The TNDP algorithm:*

TNDP algorithm uses another data structure to store the set of diagonals. With the new data structure, the above 2 operators can be implemented to run faster.

We view a diagonal as a directed one which go from the vertex with greater index to the vertex with smaller index. (Ex. *(6,4) (4,1) (3,1)*). For each vertex, we store the (increasingly) sorted list of vertices with smaller index that is connected to it by a diagonal. We have an array of lists ( Ex. *((6,[4]), (4,[1]), (3,[1]))* ).

We have to maintain this data structure for a triangulation to compute the above 2 operator.

To compute the *father(.)* operator, we have to find the *pivot* point of the edge *(n, n-1)* which is the point *p* such that *(n,a)* and *(n-1,a)* belong to the set of diagonals or an edge of the polygon. After computing the *father(.)* operator, these two diagonal or edge will be contracted to become only one *(n-1,a)*. This pivot point is exactly the last point in the list of *n* and also be the first point in the list of *n-1*. So we can compute *father(.)* in O(1) time by merely concatenate the list of *n* to the head of the list of *n-1* and delete the list of vertex *n*.

To compute the $son^i(.)$ operator, we have to split the diagonal (or edge) *(n, i)* into two ones *(n+1, i), (n, i)*. The computation of TNDP contains two phases. In first phase, we go bottom up until $P_1=P_2$. In the second phase, we go top down and compute the *son* operator. We compute the $son^i(.)$ operator by traverse from the tail of the list of vertex n to the head of it, whener we meet the vertex i in this list, we can split the list at this point. Assume the list of vertex *n* is *[a,...,i,j,...,b]*, after computing this operator the list of vertex *(n+1)* will be *[a,...,i]* and the list of vertex *n* will be *[j,...,b]*. All the diagonals corresponding to the list *[j,...,b]* will be fixed. So in one time it is called, the operator soni will fix some diagonals. But there are only *n-3=O(n)* diagonals in a triangulation. Thus, the overall time we have to traverse through the list of vertex *n* is only O(n).

*Theorem 3: TNDP* algorithm runs in time O(n).

*Proof:*

According to the above argument, we have the time of computing an operator *father(.)* is O(1), the overall time to compute $son^i(.)$ operators is *O(n)*. We have to go through at most *(n-2)* levels of the hierachy tree. Thus, the time complexity of TNDP algorithm is *O(n)*.

*Theorem 4: TNDP* algorithm is optimal.

*Proof:*

Trivially, the input for TNDP algorithm is a pair of triangulation which size is *O(n)* because each triangulation contain *(n-3)* diagonals. Thus, at least the algorithm must read all its input, it takes *O(n)* time to do this. In theorem 3, we show that TNDP can run in *O(n)* time. So this is an optimal algorithm. *Q.E.D*